\documentclass[12pt]{article}
\usepackage[utf8]{inputenc}
\usepackage[english]{babel}
\usepackage{color}
\usepackage{graphicx}
\usepackage{vmumath}
\usepackage{amsmath}
\usepackage{amsthm}
\usepackage{avo_book}
\usepackage{hyperref}
\usepackage{caption}
\usepackage{booktabs}
\usepackage{subcaption}
\usepackage[a4paper,left=15mm,right=15mm,top=30mm,bottom=20mm]{geometry}
\noaftermath
\begin{document}

\bibliographystyle{unsrt}

\title{Enumeration of $r$-regular Maps on the Torus. \\Part II: Enumeration of Unsensed Maps}

\author{
Evgeniy Krasko \qquad  Alexander Omelchenko\\
\small St. Petersburg Academic University\\
\small 8/3 Khlopina Street, St. Petersburg, 194021, Russia\\
\small\tt \{krasko.evgeniy, avo.travel\}@gmail.com
}

\begin{abstract}
The second part of the paper is devoted to enumeration of $r$-regular toroidal maps up to all homeomorphisms of the torus (unsensed maps). We describe in detail the periodic orientation reversing homeomorphisms of the torus which turn out to be representable as glide reflections. We show that considering quotients of the torus with respect to these homeomorphisms leads to maps on the Klein bottle, annulus and the M\"obius band. Using $3$- and $4$-regular maps as an example we describe the technique of enumerating quotient maps on surfaces with a boundary. Obtained recurrence relations are used to enumerate unsensed $r$-regular maps on the torus for various $r$.
\end{abstract}

\maketitle

\section*{Introduction}

\begin{figure}[ht]
\centering
	\centering
    	\includegraphics[scale=0.7]{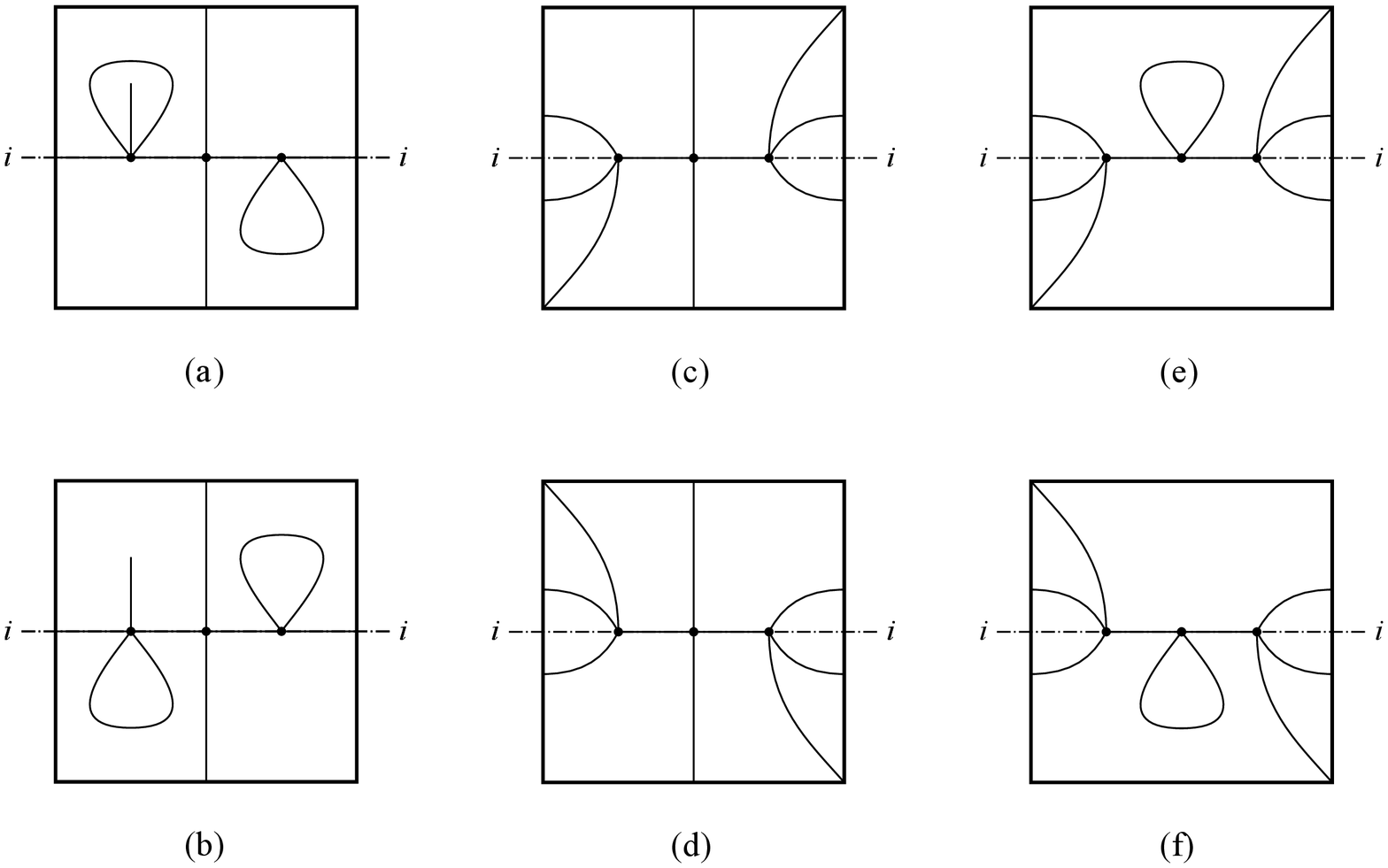}
	\caption{}
\label{fig:torus_sym_refl_examples}
\end{figure}

Following \cite{Wormald} and \cite{Walsh_generation}, we define equivalence classes of maps up to all homeomorphisms, including orientation-reserving, as unsensed maps. Taking the complete group of symmetries into account reduces the number of non-isomorphic maps. For example, the Atlas \cite{An_atlas} mentions $23$ different $4$-regular toroidal map with 3 vertices. At the same time among these maps there are three pairs which differ only by a reflection (Figure \ref{fig:torus_sym_refl_examples}(a-f)). Consequently, there are $20$ different unsensed $4$-regular maps with $3$ vertices on the torus.

Unlike the problem of counting sensed maps the problem of enumerating unsensed ones is little investigated. An exception in this regard is the paper \cite{Wormald} where an algorithmic method for enumerating planar maps up to arbitrary symmetries was developed, and the work \cite{Liskovets_reductive_technique} in which the author extended his approach to counting planar maps up to homeomor\-phisms that preserve surface orientation to the ones that reverse it. 

In the present paper we employ a geometric approach to enumerate unsensed maps based on the enumeration of rooted maps on cyclic orbifolds and on the determination of the unrooted coefficients in terms of order and orientation preserving epimorphisms from orbifold fundamental groups onto cyclic groups. These approach was firstly described in the article \cite{Azevedo} in applied to enumerated unsensed maps on surfaces regardless of genus. The authors of that article pointed out that the technique developed in they paper is suitable for enumerating unsensed maps on the surfaces of given genus. In our paper we apply this technique to the specific surface, namely, to the torus. This restriction simplifies the enumeration of order and orientation preserving epimorphisms and the determination the types of orbifolds. Then we use the results obtained in the first part of the article for counting the number of some specific quotient maps on surfaces with a boundary which arise naturally when considering symmetries that reverse surface orientation. After enumerating such quotient maps we will finally be able to enumerate $r$-regular toroidal maps up to all homeomorphisms, including orientation-reversing. For $r=3$ and $r=4$ we will provide explicit enumerating formulas. For larger $r$ numerical results will be given.

\section{Enumerating $r$-regular maps on the torus up to all homeomorphisms}

To enumerate the maps with respect to all surface homeomorphisms we will use the approach which was described  in \cite{Azevedo}. It was proved \cite[Lemma 4.2]{Azevedo} that the number $\bar{\tau}_n$ of unsensed maps with $n$ edges can be calculated by the formula
\begin{equation}
\label{eq:orient_rev_maps}
\bar{\tau}_n= \dfrac{\tilde{\tau}_n+\ddot{\tau}_n}{2},
\end{equation}
where $\tilde{\tau}_n$ is the number of sensed maps with $n$ edges and $\ddot{\tau}_n$ is the number of maps admitting an orientation-reversing automorphism. Similarly to \cite{Mednykh_Hypermaps} we can conclude that the same equality holds for the number of $r$-regular maps of a given genus, and particularly for the number $\bar{\tau}_n^{(r)}$ of $r$-regular maps on the torus. The summand $\ddot{\tau}_n^{(r)}$ in the right-hand side then can be rewritten as
$$
\ddot{\tau}_n^{(r)}=\dfrac{1}{2n}\sum_{\substack{l\mid n\\l\cdot m=n}}\sum_{O\in {\rm Orb}^-(T/\Z_{2l})}h_O(m) \cdot
\Epi_o^+(\pi_1(O), \Z_{2l}).
$$
Here $h_O(m)$ is the number of rooted quotient maps on the orbifold $O$ with $m$ semiedges, corresponding to $r$-regular maps on the torus, $O$ runs through all orientation-reversing cyclic orbifolds ${\rm Orb}^-(T/\Z_{2l})$ of the torus with period $2l$, and $\Epi_o^{+}(\pi_1(O), \Z_{2l})$ is the number of order and orientation preserving epimorphisms from fundamental group of the orbifold $O$ onto the cyclic group $\Z_{2l}$. Using results of \cite[Theorem 2.4.4]{AutomorphismGroups} we can conclude that the Riemann--Hurwitz formula for the orbifolds of the torus has the form
$$
2l\left(\chi-\sum\limits_{i=1}^r\left(1-\dfrac{1}{m_i}\right)\right)=0\qquad\qquad \Longleftrightarrow \qquad \qquad 
\chi=\sum\limits_{i=1}^r\left(1-\dfrac{1}{m_i}\right),
$$
where $\chi$ is the Euler characteristic of the orbifold $O$, $m_i$ are branch indices of the branch points. Solving this equation for orientation-reversing orbifolds we obtain the following list of all possible orbifolds:
\begin{itemize}
\item[--] an annulus without branch points;\\
\item[--] a Klein bottle without branch points;\\
\item[--] a M\"obius band without branch points;\\
\item[--] a disc with two branch points of indices $2$;\\
\item[--] a projective plane with two branch points of indices $2$.
\end{itemize}
Using the techniques given in \cite{Azevedo} we can compute the numbers $\Epi_o^{+}(\pi_1(O), \Z_{2l})$ of order and orientation preserving epimorphisms for the orbifolds listed above. These numbers are equal to:
\begin{itemize}
\item[--] $\phi(l)$ if $l$ is odd and $0$ if $l$ is even, for an annulus without branch points, $\phi$ being the Euler's function;\\
\item[--] $\phi(l)$ if $l$ is odd and $4\phi(l)$ if $l$ is even, for a Klein bottle without branch points;\\
\item[--] $\phi(l)$ if $l$ is odd and $0$, $l$ if is even, for a M\"obius band without branch points;\\
\item[--] $0$ for a disc with two branch points of indices $2$;\\
\item[--] $0$ for a projective plane with two branch points of indices $2$.
\end{itemize}

\begin{figure}[ht]
\centering
	\centering
    	\includegraphics[scale=0.7]{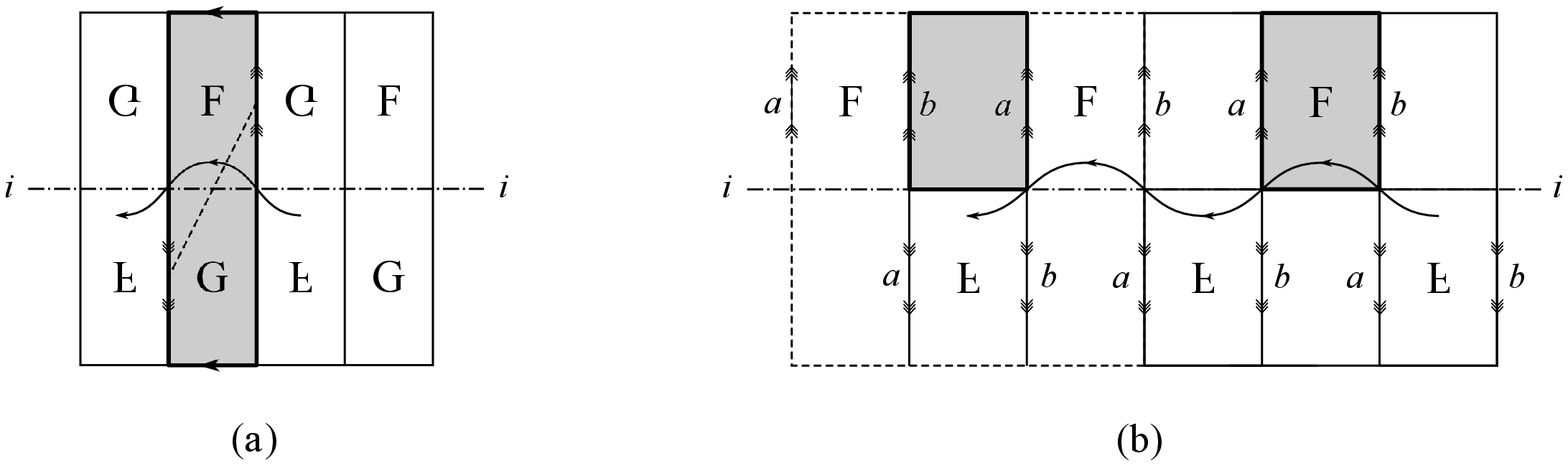}
	\caption{}
\label{fig:Glide_reflection_horis}
\end{figure}

These results admit an obvious geometric interpretation. Any orbifold corresponding to non-zero number of epimorphisms could be represented as a quotient of the torus with respect to some glide reflection of the square representing this torus on the plane. Consider first glide reflections with respect to a horizontal axis $i$ (Figure \ref{fig:Glide_reflection_horis}). Let the ratio between the value of the `shift' and the length of the side be a rational number $p/q$, $0\leq p/q<1$, $p$ and $q$ are coprime. On Figure \ref{fig:Glide_reflection_horis} an example of a glide reflection with respect to a horizontal axis and $p/q=1/4$ is shown. The fundamental polygon in this case is one fourth of a square. Since the right side of this polygon is transformed into its left side under the glide reflection with respect to $i$, we may think of these sides as glued together in the reverse direction. Consequently, the fundamental region as a whole is glued into the Klein bottle.

An example of a glide reflection with $p/q=1/3$ is shown on the Figure \ref{fig:Glide_reflection_horis} (b). For such a ratio of $p$ and $q$ the fundamental region is one sixth of the square (shaded area on Figure \ref{fig:Glide_reflection_horis}(b)). Indeed, it would take six steps for the glide reflection to transform it into itself. After the second step its left side $a$ would coincide with its right side $b$, and vice versa after the fourth step. At the same time its top and bottom sides would never become coincident. Consequently this glide reflection corresponds to a rectangular fundamental region with its right and left sides glued together. In other words, it is an annulus.

An analogous situation takes place in the general case. Namely, if $q$ is an even number then the fundamental region is a $q$-th part of the square glued into a Klein bottle (Figure \ref{fig:Glide_reflection_horis}(a)). In this case the number of possible different values for $p$ is equal to $\phi(q)$, and $q=2l$. If $q$ is odd then the fundamental region is one $2q$-th of the square glued into an annulus (Figure \ref{fig:Glide_reflection_horis}(b)). For a given $q$ the number of possible coprime values of $p$ is still equal to $\phi(q)$, and $q=l$. Conversely, for a given $l=2d$ the number $q=2l=4d$, $\phi(q)=2\phi(l)=2\phi(2d)$, and the corresponding orbifold is the Klein bottle. In the case $l=2d+1$ the number $q$ is equal to either $l$ or $2l$. In both cases $\phi(q)=\phi(l)$, but in the former case the orbifold $O$ is an annulus and in the latter case $O$ is the Klein bottle. 

\begin{figure}[ht]
\centering
	\centering
    	\includegraphics[scale=0.7]{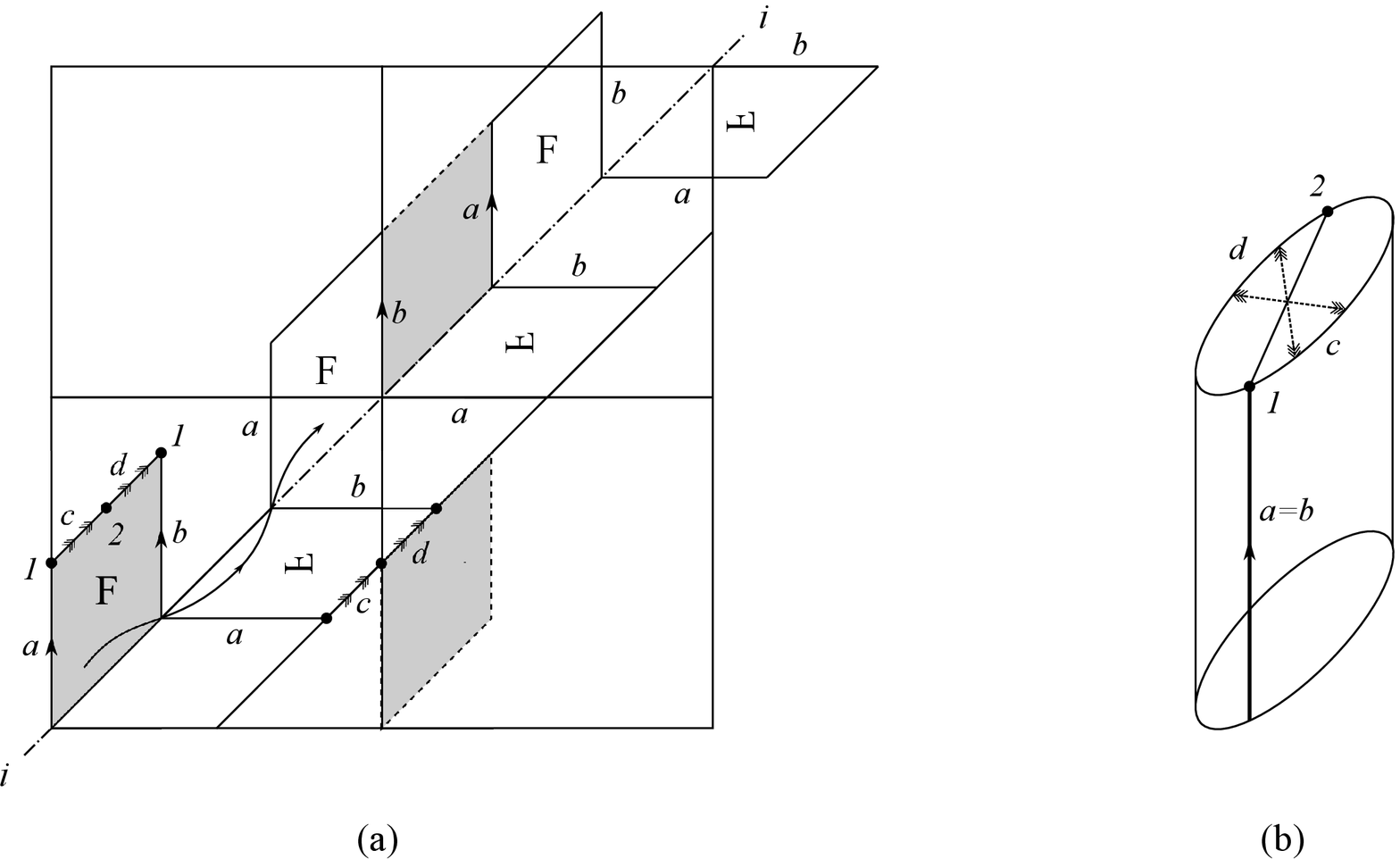}
	\caption{}
\label{fig:Glide_reflection_diag}
\end{figure}

Now consider glide reflections with respect to the diagonal of the square. Consider first the case of $q=3$, $p=1$ (Figure \ref{fig:Glide_reflection_diag},a). As in the case of a horizontal glide reflection the fundamental region is equal to one sixth of the square (shaded region $F$ one the bottom left, Figure \ref{fig:Glide_reflection_diag}(a)). In this case the sides $a$ and $b$ are glued, as well as the segments $c$ and $d$ of the side that doesn't lie on the axis. Indeed, after applying the glide reflection that moves the fundamental region in the upper right direction, its side $b$ will coincide with $a$ on the third step and vice versa (see the dashed region on the top right square). Now observe that after applying one step of the glide reflection the segment $d$ becomes coincident with the segment $c$ of the same fundamental region drawn on the bottom right square (see the dashed region on Figure \ref{fig:Glide_reflection_diag}(a)). Consequently, this glide reflection glues the fundamental region into a M\"obius band. Indeed, gluing the sides $a$ and $b$ results in an annulus (Figure \ref{fig:Glide_reflection_diag},b), and gluing the segments $c$ and $d$ means that the opposite points of one of its bases are identified. Such a surface is a M\"obius band.  

\begin{figure}[h]
\centering
\begin{tabular}[t]{cc}
	\begin{subfigure}[b]{0.5 \textwidth}
	\centering
    		\includegraphics[scale=0.7]{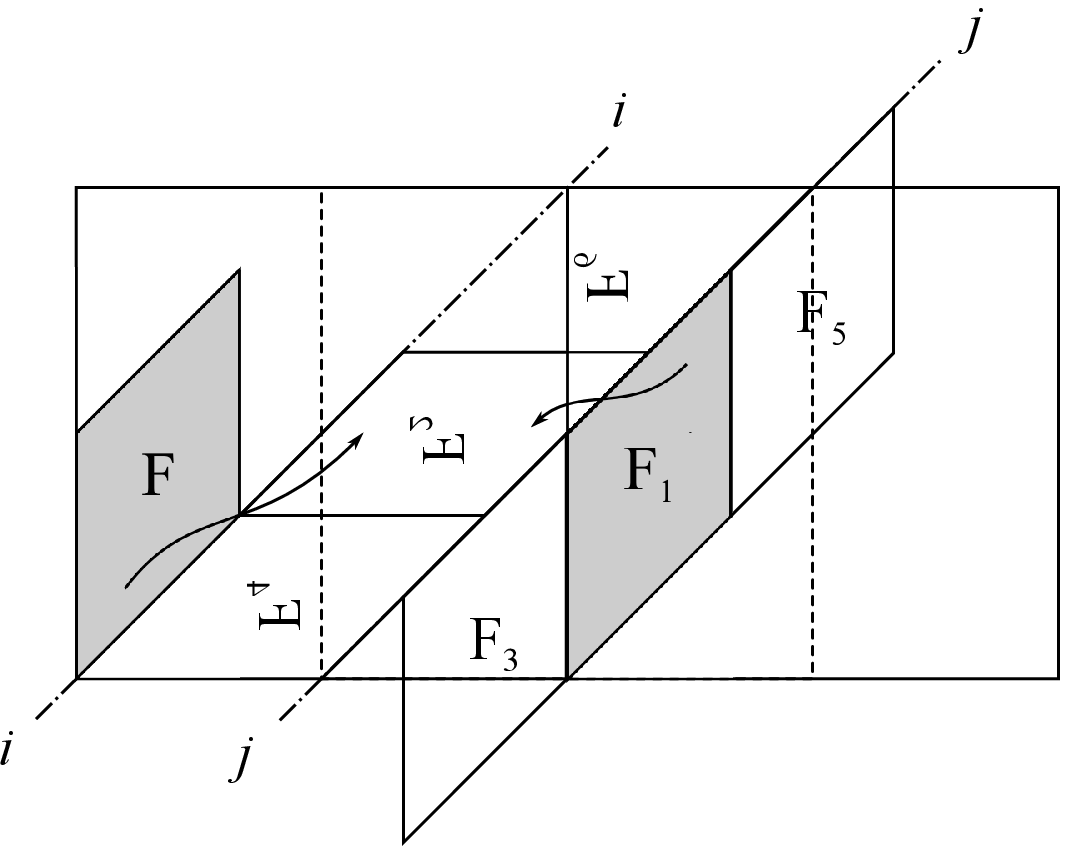}
 	\caption{}
	\end{subfigure}

	\begin{subfigure}[b]{0.5 \textwidth}
	\centering
    		\includegraphics[scale=0.7]{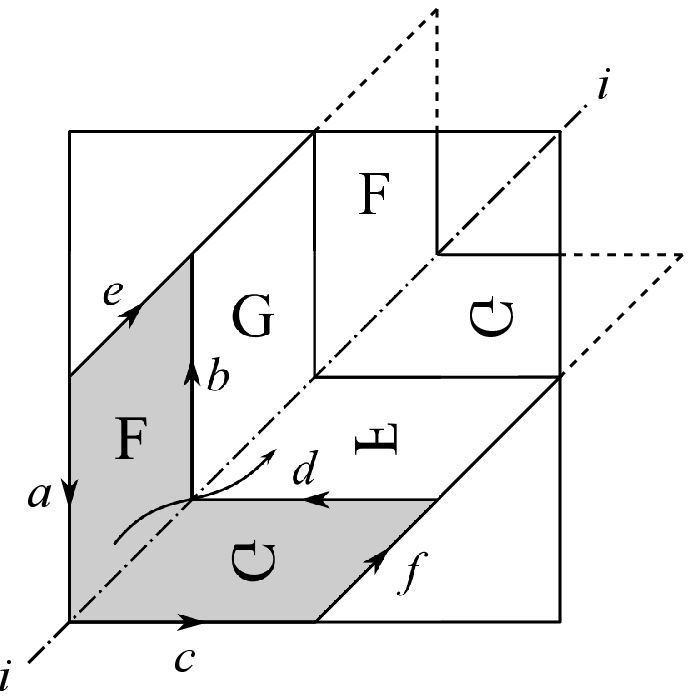}
 	\caption{}
	\end{subfigure}
\end{tabular}
\caption{}
\label{fig:Glide_reflection_diag_1}
\end{figure}

It turns out that the diagonal glide reflections corresponding to an odd value  $q$ and an even value  $2q$ actually define the same tranformation. Indeed, consider the abovementioned example of the glide reflection with respect to the axis $i$ corresponding to $q=3$ (Figure \ref{fig:Glide_reflection_diag_1}(a)). For the values $p=2$, $q=3$ after one application of the glide reflection with respect to $i$ the fundamental region $F$ is moved upwards to the region $F_2$. Note that instead of the original square we can represent the torus as another square which is obtained from it by a shift to the right by a half of the side of this square  (the dashed square on Figure \ref{fig:Glide_reflection_diag_1}(a)). Choose $F_1$ as the new fundamental region in this square and $j$ as the new glide reflection axis (Figure \ref{fig:Glide_reflection_diag_1}(a)). For the parameters $p=5$, $q=6$ of the glide reflection the fundamental region $F_1$ which is the same region as $F$ on the torus is transformed into the region $F_2$ as before (Figure \ref{fig:Glide_reflection_diag_1}(a)). Consequently, two considered glide reflections on the torus are equivalent. So for odd $l=2d+1$ we have $q=l$, $\phi(q)=\phi(2d+1)$ possible values of $p$, and the orbifold $O$ is a M\"obius band.

It remains to consider the case of $q$ being a multiple of four (see Figure \ref{fig:Glide_reflection_diag_1}(b) for the case $q=4$). The fundamental region (dashed area $F-G$ on Figure \ref{fig:Glide_reflection_diag_1}(b)) is one fourth of the square which represents the torus. The glide reflection with respect to the axis $i$ which shifts this axis by one fourth of the diagonal of the square transforms the side $c$ of the fundanetnal region into its side $b$, and the side $a$ into the side $d$ (the point $1$ is transformed into the point $2$). The next application of the glide reflection glues the sides $e$ and $f$. As a consequence, the fundamental region is glued into the Klein bottle. So in this case we have $2l=q$ and the number of possible values of $p$ is equal to $\phi(q)=\phi(2l)=2\phi(2d)$. 

As a consequence, the formula (\ref{eq:orient_rev_maps}) for the number of $r$-regular unsensed maps can be rewritten as
\begin{equation}
\label{eq:r_reg_maps_torus}
\begin{aligned}
\bar{\tau}_n^{(r)}&=
\dfrac{\tilde{\tau}_n^{(r)}}{2}+
\dfrac{1}{4n}\sum\limits_{(2d+1)\mid n}\phi(2d+1)\cdot \left[C^{(r)}_{2n/(2d+1)}+M^{(r)}_{2n/(2d+1)}\right]+\\
&+\dfrac{1}{4n}\sum\limits_{(2d+1)\mid n}\phi(2d+1)\cdot \kappa^{(r)}_{n/(2d+1)}+\dfrac{1}{n}\sum\limits_{2d\mid n}\phi(2d)\cdot \kappa^{(r)}_{n/(2d)}.
\end{aligned}
\end{equation}
Here $\kappa^{(r)}_{n}$ is the number of $r$-regular maps on the Klein bottle which was determined in the first part of the paper. $C^{(r)}_n$ and $M^{(r)}_n$ are the numbers of quotient maps with $n$ semi-edges on the annulus and on the M\"obius band.

Our next step is to enumerate the corresponding quotient maps on such orbifolds. Maps on the Klein bottle were already enumerated in the first part of the paper. The next sections are devoted to enumerating maps on surfaces with a boundary.

\section{Enumeration of maps on surfaces with a boundary}

\begin{figure}[ht]
\centering
	\centering
    	\includegraphics[scale=0.7]{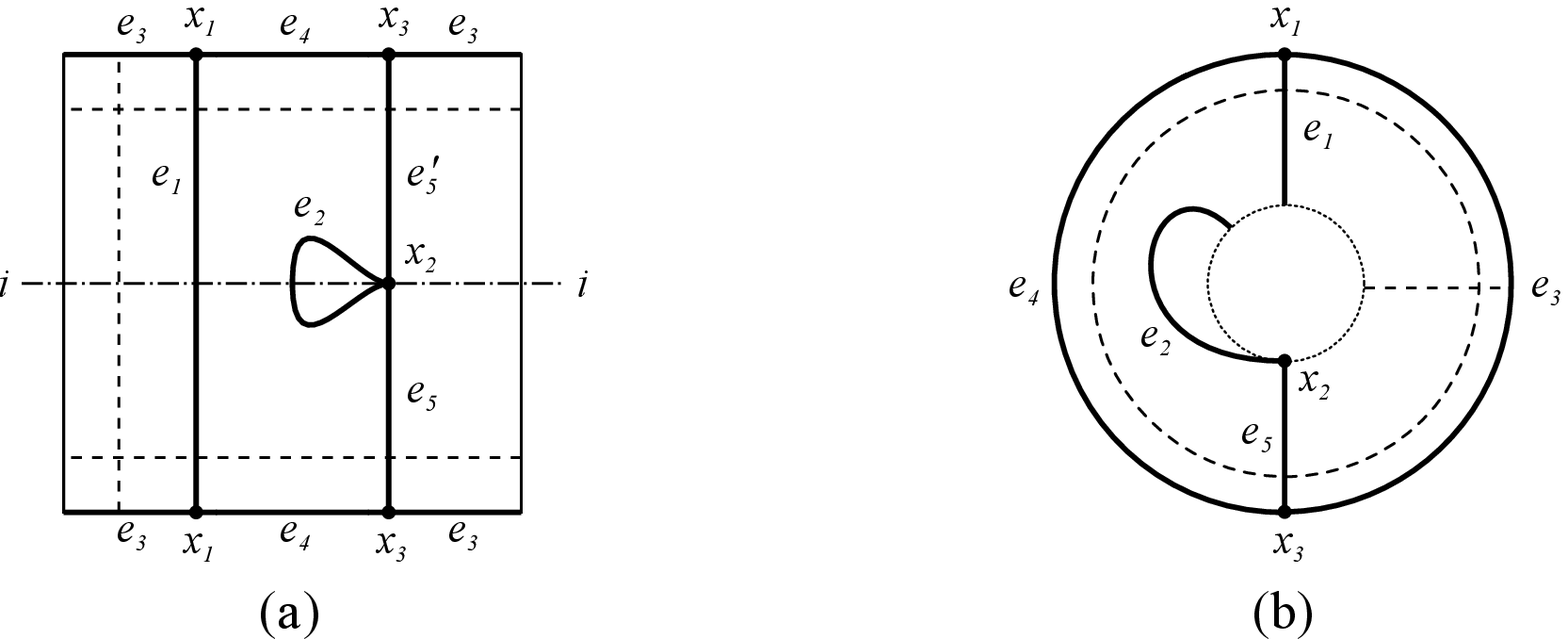}
	\caption{}
\label{fig:surfaces_with_boundaries}
\end{figure}

As in the case of homeomorphisms that  preserve surface orientation, we will enumerate quotient maps $\mathfrak{M}$ on orbifolds $O$, that is, maps that have some peculiarities that are related to the existence of the boundary of a surface. As an example consider the $4$-regular map $\M$ on the torus depicted on the Figure \ref{fig:surfaces_with_boundaries}(a). Such map is transformed into itself under the action of the reflection with respect to the axis  $i$. The corresponding quotient map $\mathfrak{M}$ on an orbifold $O$ is actually a quotient map on the annulusannulus (Figure \ref{fig:surfaces_with_boundaries}(b)). It can be seen that the boundary of an orbifold may contain vertices (vertices $x_1$, $x_2$ and $x_3$ on Figure \ref{fig:surfaces_with_boundaries}(b)) and/or edges (edges $e_3$ and $e_4$ on Figure \ref{fig:surfaces_with_boundaries}(b)). In addition, the edges $e_1$ and $e_2$ on the torus that cross the axis of symmetry correspond to dangling semi-edges of the quotient map $\mathfrak{M}$ on the annulus. Let $L$ be the period of the symmetry, for our example $L=2$. Any internal vertex (edge) of a quotient map $\mathfrak{M}$ corresponds to $L$ vertices (edges) of the map $\M$ on the torus. Any vertex (edge) that lies on the boundary of an orbifold corresponds to $L/2$ vertices (edges) of the corresponding map $\M$. Analogously, any dangling semi-edge of $\mathfrak{M}$ corresponds to $L/2$ regular edges of the map $\M$. In addition to that, the degrees of vertices that lie on the boundary of the orbifold are smaller than those of $\M$. It is convenient to think that every semi edge that lies on the boundary is counted as $1/2$ both in the number of semi-edges of the quotient map and in the degree of the vertex incident to it. With this approach the total number of semi-edges of a quotient map $\mathfrak{M}$ becomes equal to one $L$-th of the number of semi-edges of the original map $\M$ and the degree of every vertex lying on a boundary becomes equal to the one of the corresponding vertex of $\M$ divided by 2.

The next limitation follows from the fact that every face of a map on the torus must be homeomorphic to a disc. In other words we can't walk along the meridian or parallel without intersecting an edge or a vertex of the map. Now assume that a map $\M$ possesses a symmetry under which the corresponding orbifold is an annulus (Figure \ref{fig:surfaces_with_boundaries}). The requirement that forbids walking along a meridian or parallel on a torus without crossing the edges and vertices corresponds to the requirement that any walk that either wraps the annulus or passes from one of its boundaries to another must cross an edge or a vertex of the map (see Figure \ref{fig:surfaces_with_boundaries} and forbidden paths shown in dashed lines).

Now consider the situation when one edge $e$ of a quotient map $\mathfrak{M}$ lies on the boundary entirely. Note that it is only possible if there is at least one vertex $x$ lying on the same boundary. Depending on the parity of the vertex degrees of the $r$-regular map on the torus there are two different cases. For an even $r$ any vertex $x$ that lies on the boundary must be incident either to zero or to two semi-edges lying on that boundary. For an odd $r$ it must be incident to exactly one such semi-edge.

It follows that for even $r$ the existence of at least one semi-edge that lies on a given boundary means that this whole boundary is covered with edges (Figure \ref{fig:surfaces_with_boundaries}). This observation allows to simplify the enumeration of quotient maps on a surface with a boundary for even values of $r$. If some boundary of a surface is entirely covered with edges, we can contract all these edges into a single vertex with the degree uniquely determined by the number $k$ of vertices lying on this boundary. For example, for $4$-regular maps this degree is equal to exactly $k$. 

For odd values of $r$ any vertex that lies on a boundary has exactly one incident semi-edge that lies on this boundary. Consequently the cases of odd and even $r$ possess some different properties. In the remaining part of the article we will demostrate them using $4$ and $3$-regular maps. We will also present a table of enumeration results for $r=3,\ldots,6$.  

\section{Toroidal $4$-regular maps up to all homeomorphisms}

As it was shown before, to enumerate maps on the torus taking orientation-reversing homeomorphisms into account it is necessary to enumerate quotient maps on the annulus and on the M\"obius band. Since for $4$-regular maps each boundary of a surface is either fully covered by map edges or contains no map edges at all, it will be natural to classify such quotient maps with respect to the number of boundaries covered by edges. Namely, we can express the numbers $C_n$ and $M_n$ as
$$
C_n^{(4)} = C_{0,n} + C_{1,n} + C_{2,n},\qquad\qquad M_n^{(4)} = M_{0,n} + M_{1,n},
$$
$C_{i,n}$ and $M_{i,n}$ being the numbers of quotient maps on the corresponding surfaces with exactly $i$ boundaries covered by edges.

\begin{figure}[ht]
\centering
	\centering
    	\includegraphics[scale=0.7]{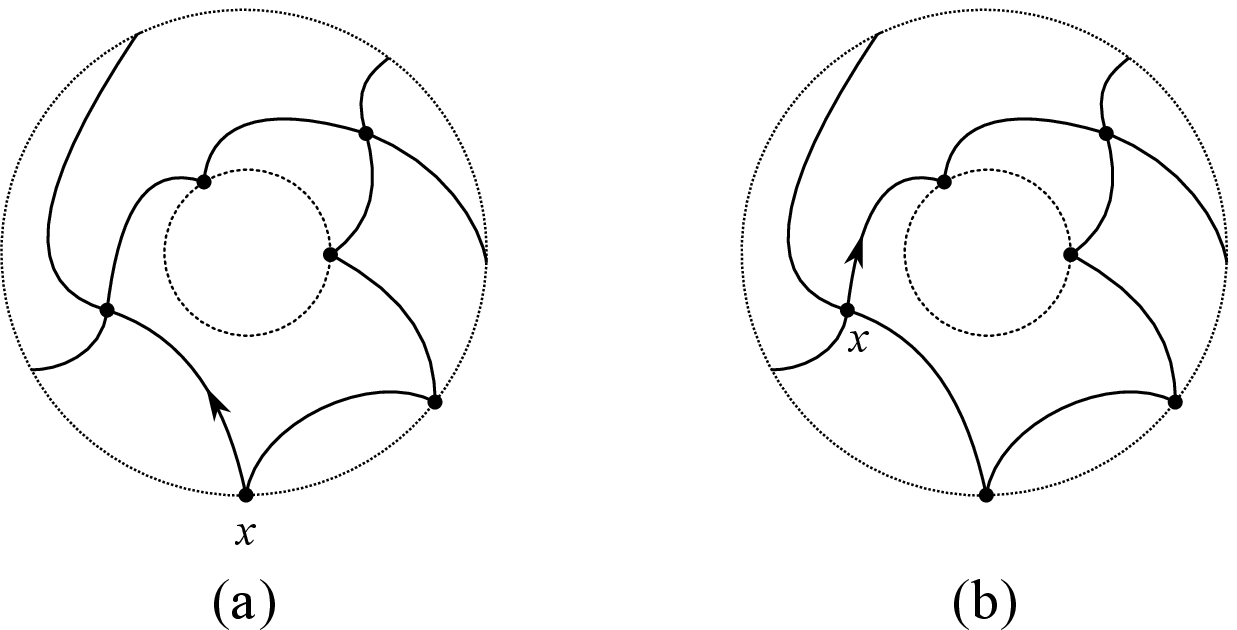}
	\caption{}
\label{fig:cylinder_0_1}
\end{figure}

Note that for a fixed $i$ the set of corresponding quotient maps on a surface with a boundary can be split further based on the position of the root vertex with respect to the boundary. We begin with the case $i=0$ (see Figure \ref{fig:cylinder_0_1}, a $4$-regular map on an annulus). For such $i$ two types of quotient maps exist: those with the root vertex lying on a boundary (Figure \ref{fig:cylinder_0_1}(a)), and those with the root vertex lying in the interior of the surface (Figure \ref{fig:cylinder_0_1}(b)). We denote by $c_{n,4}^{\rm (in)}$ and $m_{n,4}^{\rm (in)}$ the numbers of quotient maps with $n$ semi-edges on the annulus and on the M\"obius band which have the root vertex of the degree $4$ that lies in the interior of the surface (vertex $x$ on Figure \ref{fig:cylinder_0_1}(b)). By $c_{n,2}^{\rm (out)}$ and $m_{n,2}^{\rm (out)}$ we denote the numbers of quotient maps with the root vertex $x$ lying on the boundary and having the degree $2$ and the root semi-edge chosen in such a way that in the counter-clockwise ordering of the semi-edges incident to the root it is the one that appears right before the boundary (Figure \ref{fig:cylinder_0_1}(a)). Since in the general case any semi-edge could be the root, the number of quotient maps having the root vertex lying on the boundary is equal to $2\,c_{n,2}^{\rm (out)}$ and $2\,m_{n,2}^{\rm (out)}$, and the numbers $C_{0,n}$ and $M_{0,n}$ can be computed by the formulae
$$
C_{0,n} = c^{\rm (in)}_{n,4} + 2 \cdot c^{\rm (out)}_{n,2},\qquad\qquad M_{0,n} = m^{\rm (in)}_{n,4} + 2 \cdot m^{\rm (out)}_{n,2}.
$$

\begin{figure}[ht]
\centering
	\centering
    	\includegraphics[scale=0.7]{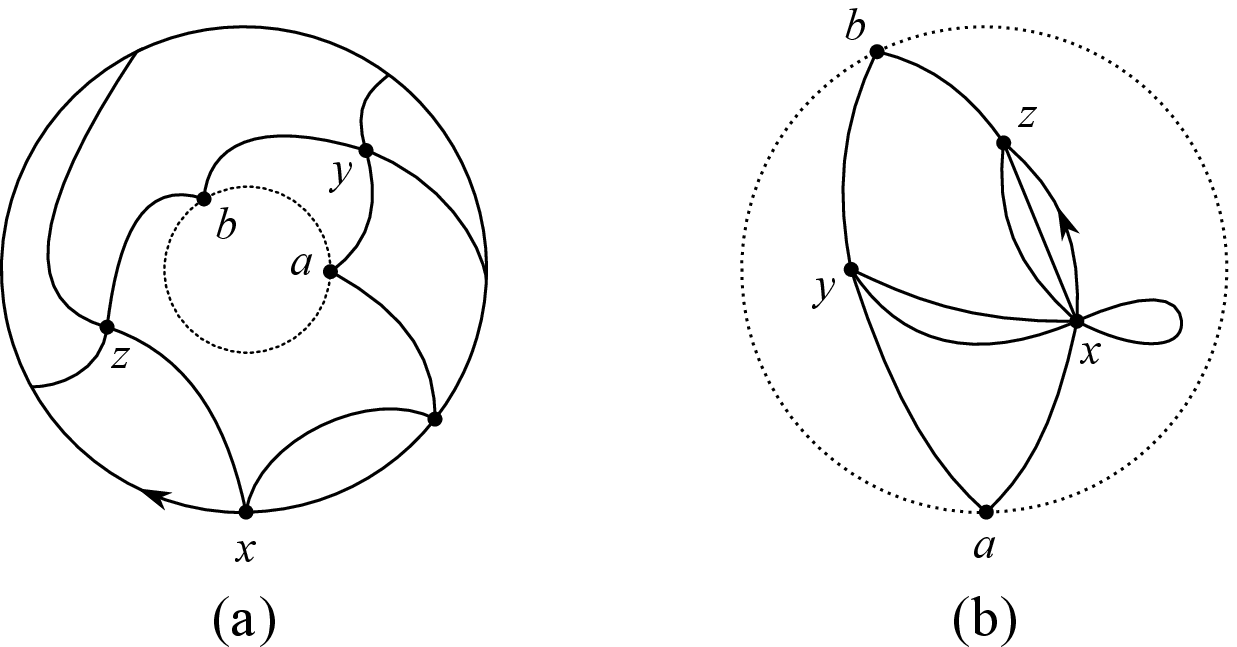}
	\caption{}
\label{fig:cylinder_1_1}
\end{figure}

Next we obtain the formulae for the numbers $C_{1,n}$ and $M_{1,n}$. Instead calculating them directly it will be more convenient to count the numbers $C_{1,n}^k$ and $M_{1,n}^k$ of quotient maps having $n$ semi-edges, exactly one boundary covered with edges, $k$ edges on this boundary and the root edge chosen in such a way that the surface is located to the right of it (see Figure \ref{fig:cylinder_0_1}(a) corresponding to the annulus, $n=27$ and $k=7$). The numbers $C_{1,n}$ and $M_{1,n}$ can be expressed through $C_{1,n}^k$ and $M_{1,n}^k$ using double counting. Indeed, for a fixed $k$ quotient maps enumerated by the sequences $C_{1,n}^k$ and $M_{1,n}^k$ have $k$ possible positions of the root semi-edge, whereas the corresponding quotient maps with an arbitrary root position have $n$ such positions. Since the parameter $k$ may take any value in the range from $1$ to $n$, we conclude that
$$
C_{1,n} = \sum_{k=1}^n \dfrac{n}{k} \cdot C_{1,n}^k,\qquad\qquad
M_{1,n} = \sum_{k=1}^n \dfrac{n}{k} \cdot M_{1,n}^k.
$$

In order to find an expression for the $C_{2,n}$ we introduce the sequence $C_{2,n}^{k, l}$ enumerating quotient maps with two distinct semi-edges located on two different boundaries of the annulus in a way that the interior of the annulus remains to the right of them. The parameters $k$ and $l$ are defined as the total numbers of semi-edges on the corresponding boundaries. Using the double counting principle we conclude that
$$
C_{2,n} = \sum_{k=1}^{n-3} \sum_{l=1}^{n-k-2}\dfrac{1}{2}\cdot \dfrac{n}{k \cdot l} \cdot C_{2,n}^{k, l}.
$$
The multiplier $1/2$ is necessary to stop distinguishing the annulus boundaries. 

As we already noted before, quotient maps that have a root edge lying on some boundary that is fully covered with edges can be reduced by contracting this whole boundary into a single vertex. As an example, the Figure \ref{fig:cylinder_1_1}(b) shows a quotient map on the disc that is obtained from the quotient map on the annulus shown on Figure \ref{fig:cylinder_1_1}(a) by contracting the edges lying on its boundary into a single vertex $x$ of degree $8$. Consequently, taking into account the change in the total amount of semi-edges, we may write the following expressions:
$$
C_{1,n}^k = \mathfrak{d}^{\rm (in)}_{n-k,k}, \qquad\qquad M_{1,n}^k = \tilde{p}_{n-k,k} 
\qquad\qquad C_{2,n}^{k, l} = \tilde{s}^{(2)}_{n-k-l,k,l}.
$$
In the first formula $\mathfrak{d}^{\rm (in)}_{n,k}$ is the number of quotient maps on the disc with no edges lying on its boundary, the root of degree $k$ and $n$ semi-edges in total.

Summarizing these considerations, we may write the following expressions for the numbers $C_n$ and $M_n$ of maps on the annulus and the M\"obius band with $n$ semi-edges:
\begin{equation}
\label{eq:cylinder_total}
C_n^{(4)}=c^{\rm (in)}_{n,4} + 2 \cdot c^{\rm (out)}_{n,2}+
\sum_{k=1}^n \dfrac{n}{k} \cdot \mathfrak{d}^{\rm (in)}_{n-k,k}+
\sum_{k=1}^{n-3} \sum_{l=1}^{n-k-2}\dfrac{1}{2}\cdot \dfrac{n}{k \cdot l} \cdot \tilde{s}^{(2)}_{n-k-l,k,l},
\end{equation}
\begin{equation}
\label{eq:moebius_total}
M_n^{(4)}=m^{\rm (in)}_{n,4} + 2 \cdot m^{\rm (out)}_{n,2}+
\sum_{k=1}^n \dfrac{n}{k} \cdot \tilde{p}_{n-k,k}.
\end{equation}
As it can be seen from the formulae (\ref{eq:cylinder_total}) and (\ref{eq:moebius_total}), to count quotient maps on the annulus and on the M\"obius band it remains to find the expressions for the numbers $c^{\rm (in)}_{n,d}$,  $c^{\rm (out)}_{n,d}$, $m^{\rm (in)}_{n,d}$,  $m^{\rm (out)}_{n,2}$ and $\mathfrak{d}^{\rm (in)}_{n,k}$. The next three subsections are dedicated to this problem.

\subsection{Enumeration of $4$-regular quotient maps on the disc}

\begin{figure}[ht]
\centering
	\centering
    	\includegraphics[scale=0.7]{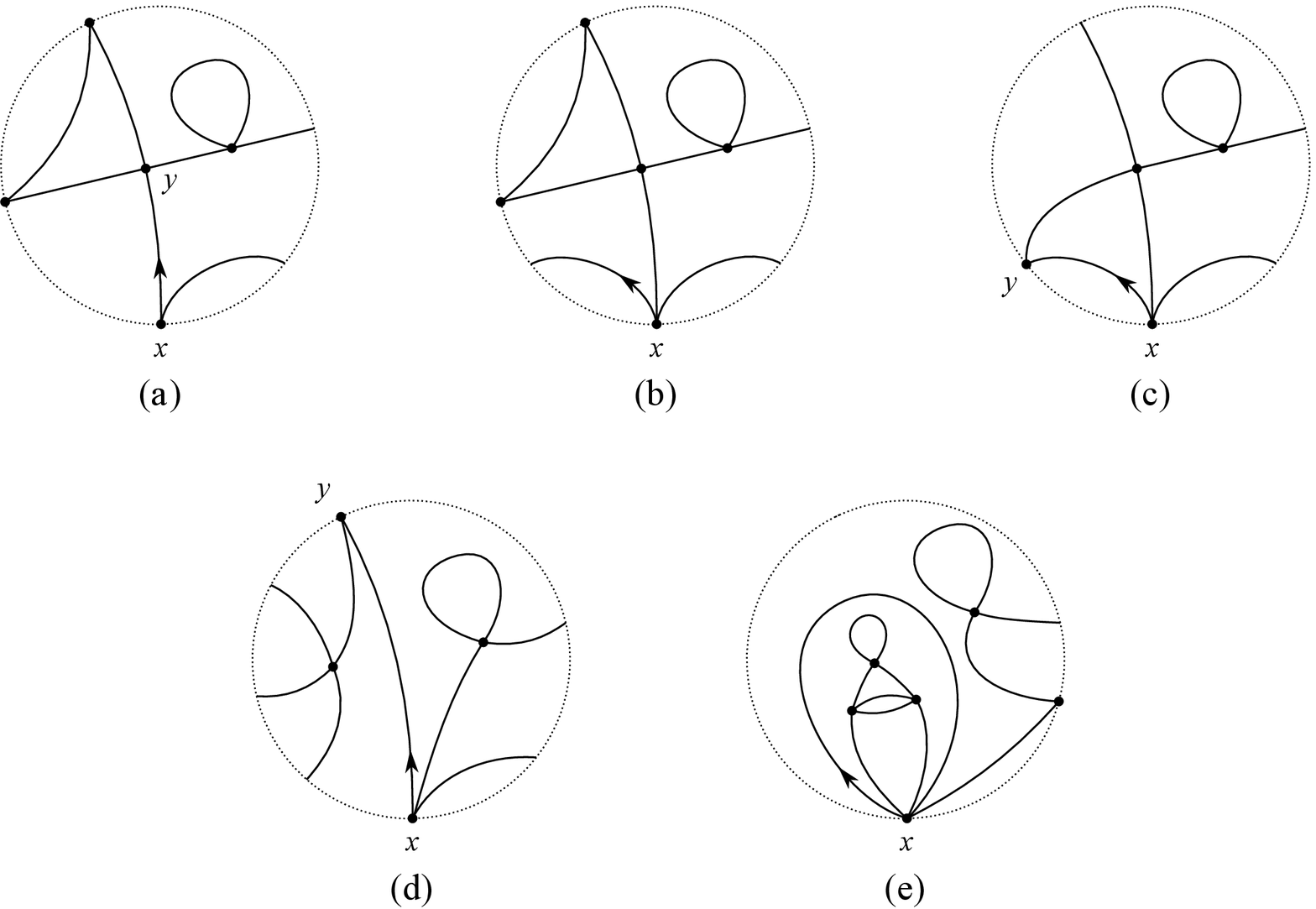}
	\caption{}
\label{fig:disk_2_1}
\end{figure}

It's the easiest to obtain recurrence relations to enumerate quotient maps on the disc. In the following sections we will need the numbers of quotient maps with the root vertex lying in the interior of the disc as well as quotient maps with the root vertex lying on its boundary. We will start with the second type, namely with the quotient maps with the root vertex but no edges lying on the boundary (Figure \ref{fig:disk_2_1}). For definiteness we will assume that the root semi-edge is the last one among semi-edges incident to the root vertex $x$ in the counter-clockwise order. Denote the number of such quotient maps by $\mathfrak{d}^{\rm (out)}_{n,d}$. To show that the recurrence relation for these numbers has the form
\begin{equation}
\label{eq:disk_2}
\mathfrak{d}^{\rm (out)}_{n,d} = \mathfrak{d}^{\rm (out)}_{n-2,d+2} + \mathfrak{d}^{\rm (out)}_{n-1,d-1} + \mathfrak{d}^{\rm (out)}_{n-2,d} + \sum\limits_{i=0}^{n-2} \mathfrak{d}^{\rm (out)}_{i,1}\cdot \mathfrak{d}^{\rm (out)}_{n-i-2,d-1} + \sum\limits_{i=0}^{n-2}\sum\limits_{j=0}^{d-2} \tilde{s}_{i,j}\cdot \mathfrak{d}^{\rm (out)}_{n-2-i,d-2-j},
\end{equation}
we analyze possible results of contracting the root semi-edge in a quotient map $\mathfrak{M}$. Let this edge join a vertex $x$ of degree $d$ with a vertex $y$ lying in the interior of the disc (Figure \ref{fig:disk_2_1}(a)). In this case contracting this edge yields a quotient map of the same type, but with the degree of the root vertex equal to $d+2$. The number of semi-edges in the obtained quotient map is reduced by two relatively to the original quotient map. This case corresponds to the summand $\mathfrak{d}^{\rm (out)}_{n-2,d+2}$ in the right hand side of (\ref{eq:disk_2}). 

Now assume that the root semi-edge is dangling and ends on the boundary of the disc (Figure \ref{fig:disk_2_1}(b)). Contracting such a semi-edge reduces by one both the total number of semi-edges and the degree of the root vertex. This case corresponds to the second summand in the right hand side of (\ref{eq:disk_2}). 

The third and the fourth summands in the right hand side of (\ref{eq:disk_2}) describe the cases when the root semi-edge joins the root vertex $x$ of the degree $d$ with a vertex $y$ of the degree $2$ lying on the boundary (Figure \ref{fig:disk_2_1}(c-d)). The second edge incident to $y$ may lie on any side of the root edge. If it's the left side (Figure \ref{fig:disk_2_1}(c)), contracting the root edge yields a quotient map with the unchanged degree of the root vertex, but the total number of semi-edges decreased by two (summand $\mathfrak{d}^{\rm (out)}_{n-2,d}$ in (\ref{eq:disk_2})). If it's the right side, then the root semi-edge splits the disc into two regions with some quotient map in each of them (Figure \ref{fig:disk_2_1}(d)). 

Finally, the last summand in the right hand side of (\ref{eq:disk_2}) corresponds to the case of the root semi-edge being a loop (Figure \ref{fig:disk_2_1}(e)). Contracting this loop separates the sub-map lying in the interior of this edge and transforms it into a map with $i$ edges and the degree of its root-vertex equal to $j$ on the sphere. The remaining part of the quotient map $\mathfrak{M}$ is transformed into a quotient map on the disc with the parameters $n-2-i,d-2-j$ (fifth summand in (\ref{eq:disk_2})). 

\begin{figure}[ht]
\centering
	\centering
    	\includegraphics[scale=0.7]{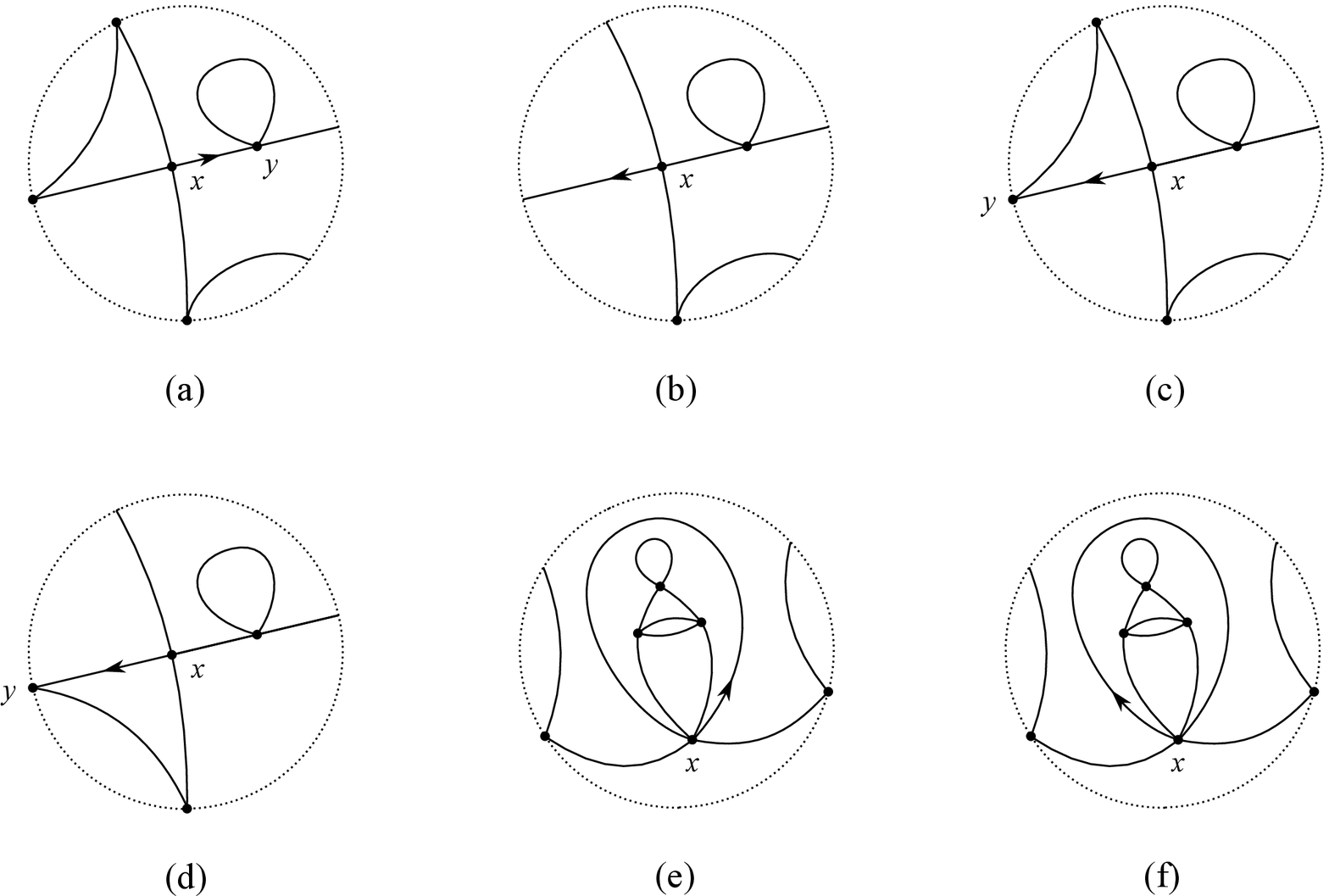}
	\caption{}
\label{fig:disk_3_1}
\end{figure}

The next step is to derive a recurrence relation for the numbers $\mathfrak{d}^{\rm (in)}_{n,d}$ of quotient maps with the root edge lying in the interior of the disc (Figure \ref{fig:disk_3_1}) and no semi-edges lying on its boundary. The numbers $\mathfrak{d}^{\rm (in)}_{n,d}$ satisfy the relation
\begin{equation}
\label{eq:disk_3}
\mathfrak{d}^{\rm (in)}_{n,d} = \mathfrak{d}^{\rm (in)}_{n-2,d+2} + \mathfrak{d}^{\rm (out)}_{n-1,d-1} + 2 \cdot \mathfrak{d}^{\rm (out)}_{n-2,d} +
2\sum\limits_{i=0}^{n-2}\sum\limits_{j=0}^{d-2} \tilde{s}_{i,j}\cdot \mathfrak{d}^{\rm (in)}_{n-2-i,d-2-j}.
\end{equation}

Indeed, the first summand in the right hand side of (\ref{eq:disk_3}) corresponds to the case of the root semi-edge joining the root vertex $x$ with a vertex $y$ lying in the interior of the disc (Figure \ref{fig:disk_3_1}(a)). The second summand corresponds to a dangling semi-edge (Figure \ref{fig:disk_3_1}(b)), contracting which moves the root vertex $x$ to the boundary of the disc. The summand $2 \cdot \mathfrak{d}^{\rm (out)}_{n-2,d}$ describes the case when the root vertex $x$ is joined by the root semi-edge with a vertex $y$ of degree $2$ lying on the disc boundary (Figure \ref{fig:disk_3_1}(c-d)). The multiplier $2$ is explained by the fact that the second edge incident to $y$ may lie on either side of the root edge. Finally, the last summand in the right hand side of (\ref{fig:disk_3_1}) corresponds to the case of the root semi-edge being a loop (Figure \ref{fig:disk_3_1}(e-f)). Since this loop may have two possible orientations, this summand in (\ref{eq:disk_3})  is multiplied by $2$.

\begin{figure}[ht]
\centering
	\centering
    	\includegraphics[scale=0.7]{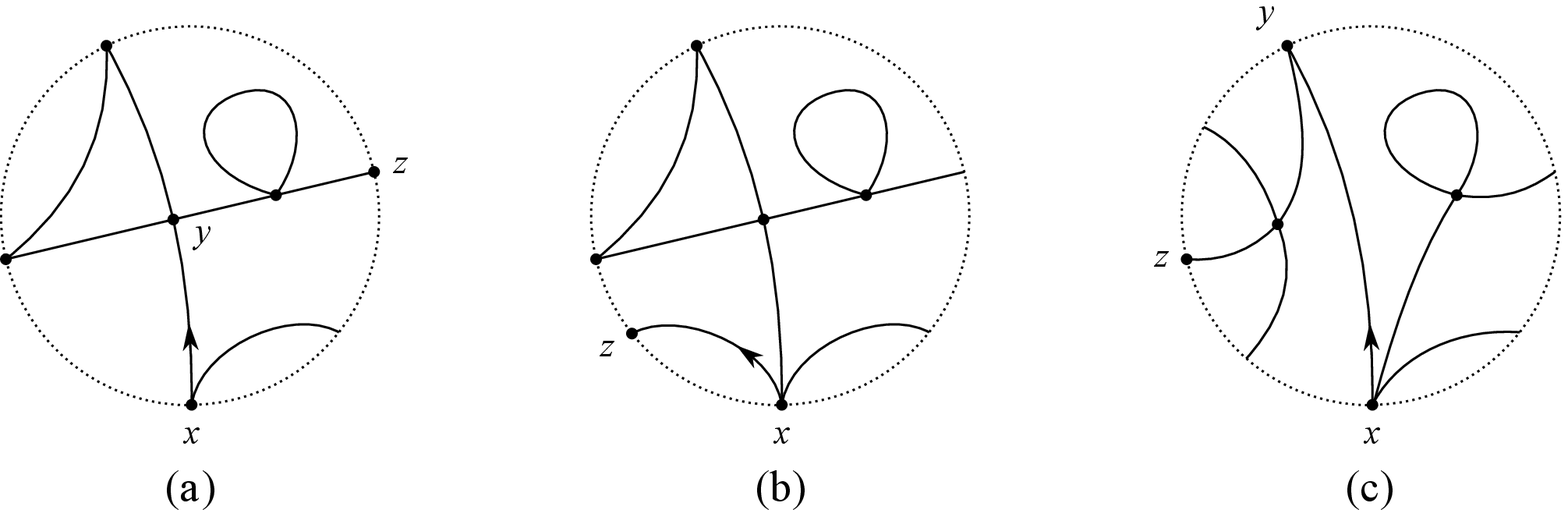}
	\caption{}
\label{fig:disk_tilde_1}
\end{figure}

To enumerate quotient maps on the M\"obius band we will also need the numbers $\mathfrak{d}^{\rm (1)}_{n,d}$ enumerating quotient maps on the disc with the root vertex $x$ lying on the boundary, no edges lying on the boundary, and one distinguished leaf $z$ lying on it (Figure \ref{fig:disk_tilde_1}(a)). The corresponding recurrence relation has the form
\begin{multline}
\label{eq:disk_tilde}
\mathfrak{d}^{\rm (1)}_{n,d} = \mathfrak{d}^{\rm (1)}_{n-2,d+2} + \mathfrak{d}^{\rm (1)}_{n-1,d-1} + 
\mathfrak{d}^{\rm (1)}_{n-2,d} +
\mathfrak{d}^{\rm (out)}_{n-2,d-1} + \sum\limits_{i=0}^{n-2} \bigl ( \mathfrak{d}^{\rm (out)}_{i,1}\cdot \mathfrak{d}^{\rm (1)}_{n-i-2,d-1}  + \mathfrak{d}^{\rm (1)}_{i,1}\cdot \mathfrak{d}^{\rm (out)}_{n-i-2,d-1} \bigr ) + \\
+\sum\limits_{i=0}^{n-2}\sum\limits_{j=0}^{d-2}\tilde{s}_{i,j}\cdot \mathfrak{d}^{\rm (1)}_{n-2-i,d-2-j}. 
\end{multline}

As compared to the relation (\ref{eq:disk_2}), a new summand $\mathfrak{d}^{\rm (out)}_{n-2,d-1}$  appears here. It corresponds to the case of the root edge joining the root vertex with the leaf $z$ (Figure \ref{fig:disk_tilde_1}(b)). In addition to that, the third summand in (\ref{eq:disk_2}) is now split into two, since we should take into account that after cutting the disc along the root edge, the leaf $z$ may remain in either of the regions (Figure \ref{fig:disk_tilde_1}(c) shows one of these possibilities).

\subsection{Enumeration of $4$-regular quotient maps on the annulus}

\begin{figure}[ht]
\centering
	\centering
    	\includegraphics[scale=0.7]{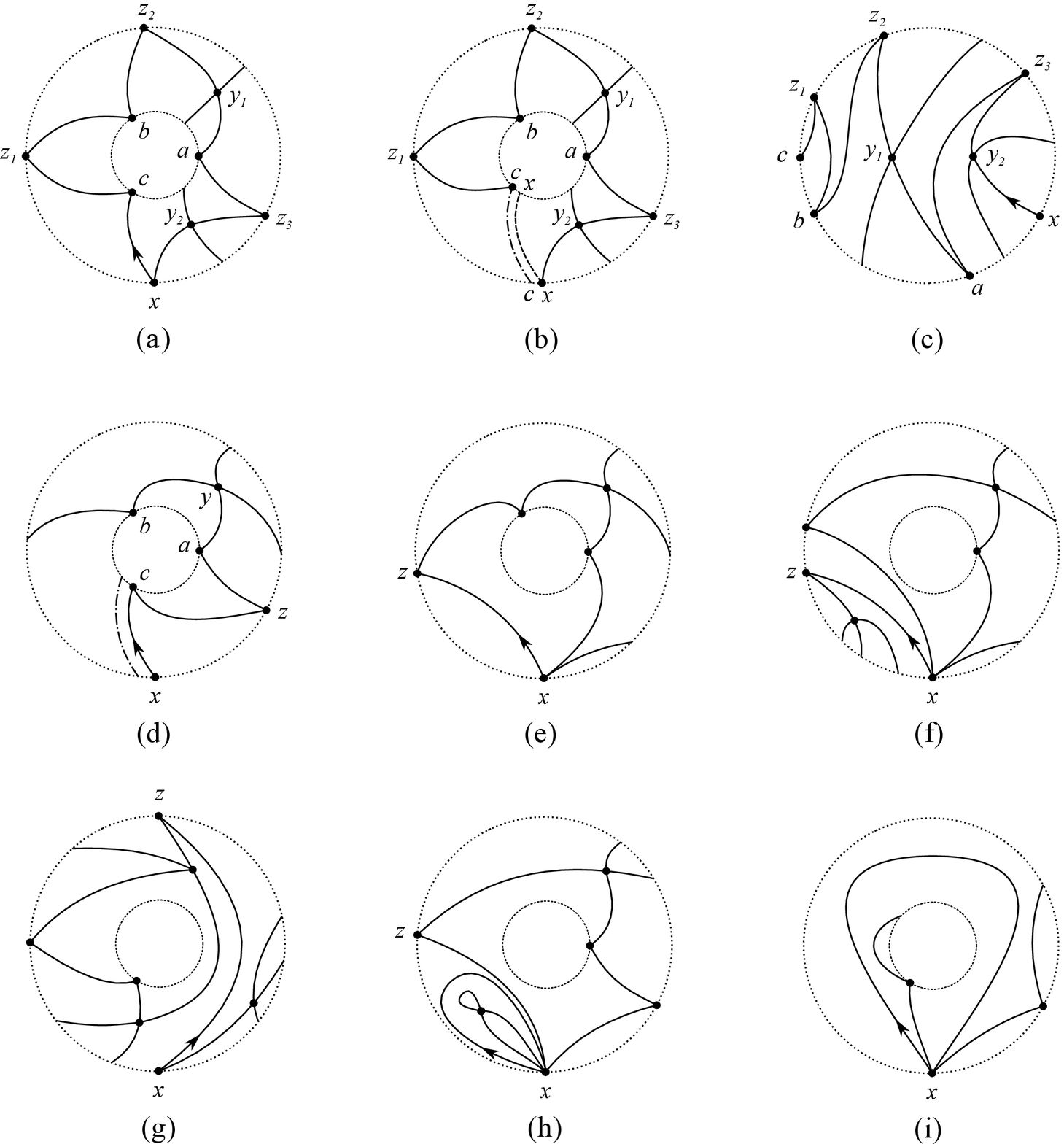}
	\caption{}
\label{fig:cylinder_2_1}
\end{figure}

Next we will derive recurrence relations for the numbers $c^{\rm (out)}_{n,d}$ and $c^{\rm (in)}_{n,d}$. We begin with the quotient maps with a root vertex lying on some boundary with no incident edges going along this boundary (Figure \ref{fig:cylinder_2_1}(a)). The number $c^{\rm (out)}_{n,d}$ of such quotient maps can be calculated by the formula
\begin{multline}
\label{eq:cylinder_2}
c^{\rm (out)}_{n,d} = c^{\rm (out)}_{n-2,d+2} + c^{\rm (out)}_{n-1,d-1}  + \mathfrak{d}^{\rm (1)}_{n-2,d-1}+ c^{\rm (out)}_{n-2,d} + 
\sum\limits_{i=0}^{n-2} \bigl ( \mathfrak{d}^{\rm (out)}_{i,1}\cdot c^{\rm (out)}_{n-i-2,d-1} + 
c^{\rm (out)}_{i,1}\cdot \mathfrak{d}^{\rm (out)}_{n-i-2,d-1} \bigr ) + \\
+\sum\limits_{i=0}^{n-2}\sum\limits_{j=0}^{d-2} \bigl (\tilde{s}_{i,j}\cdot c^{\rm (out)}_{n-2-i,d-2-j} + 
\mathfrak{d}^{\rm (out)}_{i,j}\cdot \mathfrak{d}^{\rm (in)}_{n-2-i,d-2-j} \bigr ).
\end{multline}

The combinatorial sense of the summands $c^{\rm (out)}_{n-2,d+2}$, $c^{\rm (out)}_{n-1,d-1}$ in the right hand side of (\ref{eq:cylinder_2}) is exactly the same as the one of the corresponding summands in the formula (\ref{eq:disk_2}) (see Figure \ref{fig:disk_2_1}(a-b)). The summand $\mathfrak{d}^{\rm (1)}_{n-2,d-1}$ corresponds to quotient maps with the root edge joining the root vertex with a vertex on the opposite boundary (Figure \ref{fig:cylinder_2_1}(a)). Indeed, contracting such an edge $(x,c)$ can be described as follows (Figure \ref{fig:cylinder_2_1}(b)): we split the edge $(x,c)$ into two identical edges, contract one of them into $x$ and the other into $c$. As a result the inner and outer boundaries get merged by $(x,a,b,c)$ and we obtain the quotient map shown on the Figure \ref{fig:cylinder_2_1}(c). Note that the other edge incident to $c$ in the original quotient map must lie to the right of the root edge (Figure \ref{fig:cylinder_2_1}(a)). Indeed, since the root edge is the leftmost one with respect to the root vertex $x$, if the other edge incident to $c$ was its leftmost edge (Figure \ref{fig:cylinder_2_1}(d)), on the disc there would be a path from one boundary to another that does not cross any edges (dot-dashed line on Figure \ref{fig:cylinder_2_1}(d)). But as we noted before, such an edge means that the torus contains a face not homeomorphic to a disk.

The next three summands correspond to the case of the root edge joining $x$ with a vertex $z$ of degree $2$ placed on the same boundary as $x$. In the first sub-case the other edge incident to $z$ is on the left of the root edge (Figure \ref{fig:cylinder_2_1}(e)). Contracting the root edge yields a quotient map on the annulus with $n-2$ semi-edges (summand $c^{\rm (out)}_{n-2,d}$ in (\ref{eq:cylinder_2})). The next two sub-cases correspond to this edge lying on the right of the root edge (Figure .\ref{fig:cylinder_2_1}(f-g)). This root edge, in its turn, can leave the second boundary either on its left (Figure \ref{fig:cylinder_2_1}(f)), or on its right (Figure \ref{fig:cylinder_2_1}(g)). 

The last two summands in (\ref{eq:cylinder_2}) with the summation signs correspond to quotient maps with the root edge being a loop (Figure \ref{fig:cylinder_2_1}(h-i)). The first summand enumerates maps with the root edge contractible into a point on the annulus surface (Figure \ref{fig:cylinder_2_1}(h)). The second summand enumerates the maps with the root edge wrapping the annulus boundary (Figure \ref{fig:cylinder_2_1}(i)). 

\begin{figure}[ht]
\centering
	\centering
    	\includegraphics[scale=0.7]{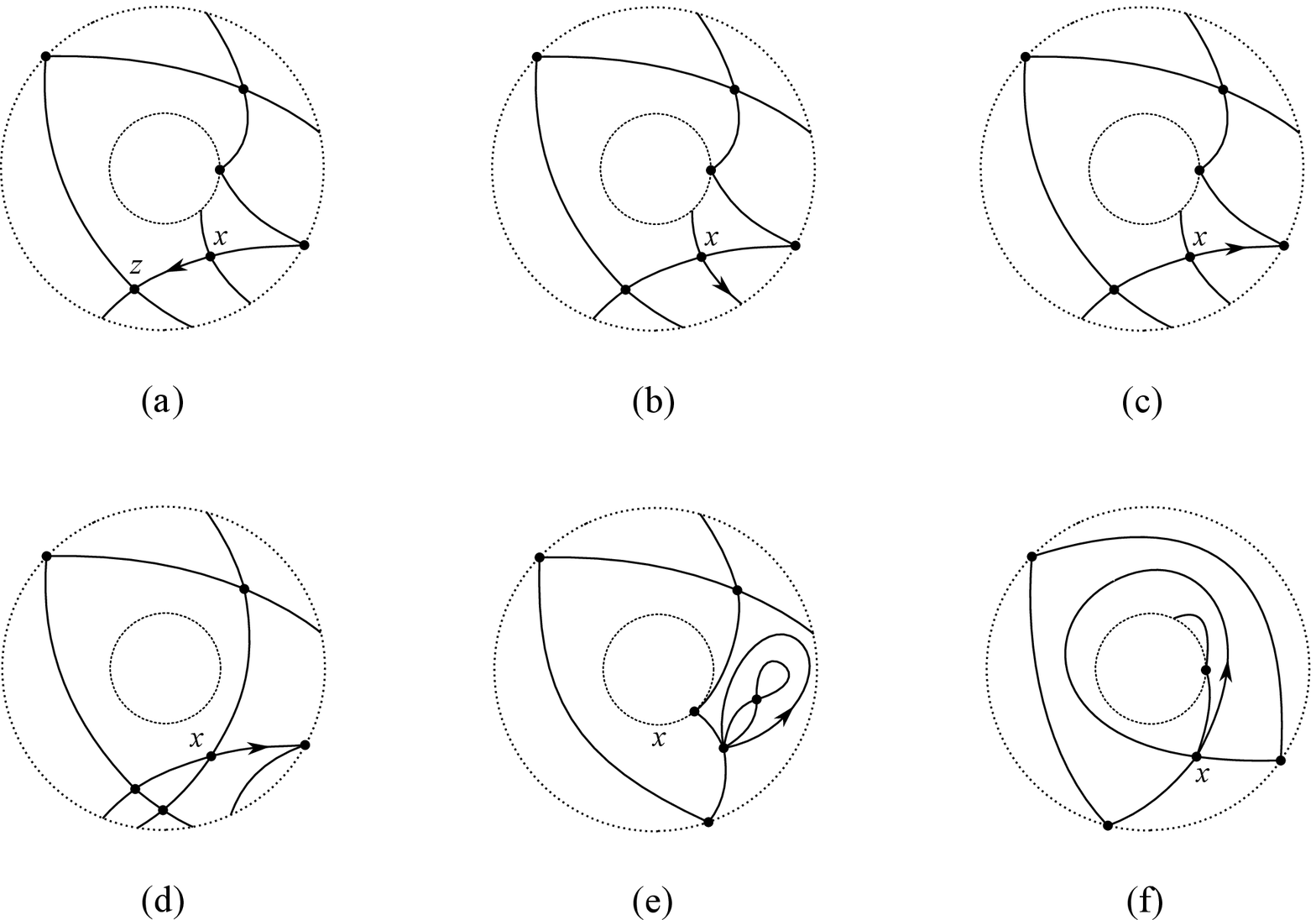}
	\caption{}
\label{fig:cylinder_3_1}
\end{figure}

Next we describe quotient maps with the root lying in the interior of the annulus. The number of such maps $c_{n,d}$ can be calculated by the the following recurrence relation:
\begin{equation}
\label{eq:cylinder_3}
c^{\rm (in)}_{n,d} = c^{\rm (in)}_{n-2,d+2} + c^{\rm (out)}_{n-1,d-1} + 2 \cdot c^{\rm (out)}_{n-2,d} + 
\sum\limits_{i=0}^{n-2}\sum\limits_{j=0}^{d-2} \bigl (2 \cdot \tilde{s}_{i,j}\cdot c^{\rm (in)}_{n-2-i,d-2-j} + d_{i,j}\cdot d_{n-2-i,d-2-j} \bigr ).
\end{equation}
The first summand in the right hand side of (\ref{eq:cylinder_3}) corresponds to the case of the root edge connecting the root vertex $x$ with one of the vertices also lying in the interior of the annulus (Figure \ref{fig:cylinder_3_1}(a)), whereas the second summand describes quotient maps with the root edge being dangling (Figure \ref{fig:cylinder_3_1}(b)). The summand $2 \cdot c^{\rm (out)}_{n-2,d}$ describes maps with the root edge connecting the root vertex $x$ with a vertex $z$ of the degree $2$ lying on the boundary of the annulus (Figure \ref{fig:cylinder_3_1}(c-d)). The multiplier $2$ is needed to take into account the fact that the second edge incident to $z$ can be positioned in two different ways with respect to the root edge (Figure \ref{fig:cylinder_3_1}(c-d)). Two summands under two summation signs enumerate quotient maps with the root edge being a loop. The first of them corresponds to the case of the root edge being contractible into a point on the surface of the annulus (Figure \ref{fig:cylinder_3_1}(e)), the multiplier $2$ reflects the fact that there are two possible orientations for such an edge. The second summand corresponds to a loop that wraps the annulus along its boundaries (Figure \ref{fig:cylinder_3_1}(f)).

\subsection{Enumeration of $4$-regular quotient maps on the M\"obius band}

\begin{figure}[ht]
\centering
	\centering
    	\includegraphics[scale=0.7]{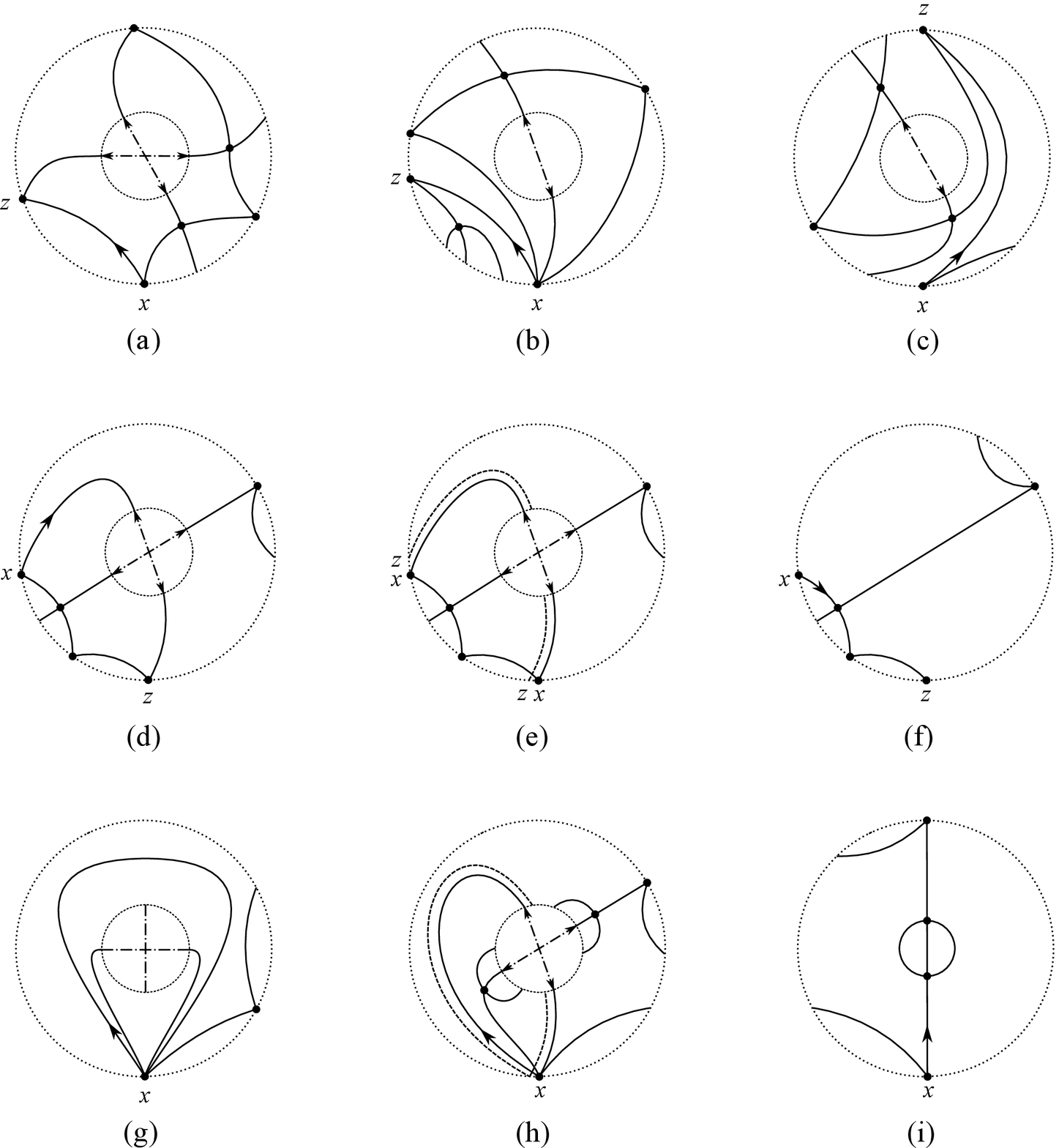}
	\caption{}
\label{fig:moebius_2_1}
\end{figure}

As in the case of the annulus we should consider two cases, one of the root vertex lying on the boundary (Figure \ref{fig:moebius_2_1}) and the other of the root vertex lying in the interior of the M\"obius band. In both cases it will be convenient for us to use the representation of the M\"obius band as an annulus with opposite points of one of its boundaries identified pairwise (Figure \ref{fig:moebius_2_1}). 

We begin with the first case. The corresponding recurrence formula for the numbers $m^{\rm (out)}_{n,d}$ takes the form
\begin{multline}
\label{eq:moebius_2}
m^{\rm (out)}_{n,d} = m^{\rm (out)}_{n-2,d+2} + m^{\rm (out)}_{n-1,d-1} + m^{\rm (out)}_{n-2,d} + 
\sum\limits_{i=0}^{n-2} \bigl ( \mathfrak{d}^{\rm (out)}_{i,1}\cdot m^{\rm (out)}_{n-i-2,d-1} + 
m^{\rm (out)}_{i,1}\cdot \mathfrak{d}^{\rm (out)}_{n-i-2,d-1} \bigr ) + \mathfrak{d}^{\rm (1)}_{n-2,d-1}+ \\
+\sum\limits_{i=0}^{n-2}\sum\limits_{j=0}^{d-2} \bigl (\tilde{s}_{i,j}\cdot m^{\rm (out)}_{n-2-i,d-2-j} + 
\mathfrak{d}^{\rm (out)}_{i,j}\cdot p_{n-2-i,d-2-j} \bigr )+ (d-1)\cdot\mathfrak{d}^{\rm (out)}_{n-2,d-2} .
\end{multline}
This formula is very similar to the analogous one (\ref{eq:cylinder_2}) for quotient maps on an annulus. In particular, the first summand in the right hand side of (\ref{eq:moebius_2}) also corresponds to contracting an edge that joins two vertices lying in the interior, whereas the second one enumerates maps with the root semi-edge being dangling (ending on the boundary of the M\"obius band).

The next four summands correspond to the case of the root edge connecting the root vertex $x$ with some vertex $z$ of degree $2$ lying on the boundary (Figure \ref{fig:moebius_2_1},(a-e)). The summand $m^{\rm (out)}_{n-2,d}$ describes maps with the second edge incident to $z$ lying on the left of the root edge (Figure \ref{fig:moebius_2_1}(a)). The next two summands, correspondingly, describe maps with this edge lying on the right of the root edge. The root edge, in its turn, may lie either on the left of the disc with pairwise identified points (Figure \ref{fig:moebius_2_1}(b)) or on the right of it (Figure \ref{fig:moebius_2_1}(c)). Finally, the root edge may also cross this disc (Figure \ref{fig:moebius_2_1}(d)). After cutting the M\"obius band along this edge (Figure \ref{fig:moebius_2_1}(e)), inverting the left part and gluing the parts back, we obtain a quotient map on the disc with an additional distinguished leaf $z$ (Figure \ref{fig:moebius_2_1}(f)).

The last three summands in the right hand side of (\ref{eq:moebius_2}) describe quotient maps with the root edge that is a loop. The first of them corresponds to maps with the root edges being loops contractible into points on the surface. By contracting such a loop we obtain a pair of quotient maps, one on the M\"obius band and another of the sphere. Maps described by the second summand have a root edge that wraps the central disc (Figure \ref{fig:moebius_2_1}(g)). Contracting such an edge again yields a pair of quotient maps, one the projective plane and another on the disc. Finally, the last summand in the right hand side of (\ref{eq:moebius_2}) describes maps in which the root edge crosses the central disc (Figure \ref{fig:moebius_2_1},h). In this case cutting along this edge, inverting one of the obtained parts and gluing these parts back yields a quotient map on a disc (Figure \ref{fig:moebius_2_1}(i)). The multiplier $(d-1)$ has the same combinatorial sense as the same multiplier in formulas that describe maps on the projective plane and on the Klein bottle.

Finally, for the numbers $m^{\rm (in)}_{n,d}$ of quotient maps on the M\"obius band with the root vertex lying in its interior the following recurrence relation holds:
\begin{multline}
\label{eq:moebius_3}
m^{\rm (in)}_{n,d} = m^{\rm (in)}_{n-2,d+2}  + m^{\rm (out)}_{n-1,d-1} + 2 \cdot m^{\rm (out)}_{n-2,d} + \\
+ 2 \sum\limits_{i=0}^{n-2}\sum\limits_{j=0}^{d-2} \bigl (\tilde{s}_{i,j}\cdot m^{\rm (in)}_{n-2-i,d-2-j} + \mathfrak{d}_{i,j}\cdot p_{n-2-i,d-2-j} \bigr )+ (d-1) \cdot  \mathfrak{d}_{n-2,d-2}.
\end{multline}
All five summands have the same combinatorial sense as the corresponding summands in the previous formula. The only difference is the presence of the multiplier $2$ for the third, fourth and fifth summands. It is needed to correctly count two possible orientations of the root loop.

\subsection{The formula for enumerating $4$-regular unlabelled maps on the torus}

Now we can use the results obtained above to enumerate $4$-regular maps on the torus up to all isomorphisms. Substituting the numbers $\mathfrak{d}^{\rm (in)}_{n,k}$, $c^{\rm (out)}_{n,d}$, $c^{\rm (in)}_{n,d}$, $m^{\rm (out)}_{n,2}$ and $m^{\rm (in)}_{n,d}$ which are calculated by the formulas (\ref{eq:disk_3}), (\ref{eq:cylinder_2}), (\ref{eq:cylinder_3}), (\ref{eq:moebius_2}) and (\ref{eq:moebius_3}) into the expressions (\ref{eq:cylinder_total}) and (\ref{eq:moebius_total}) we obtain the numbers $C^{(4)}_n$ and $M^{(4)}_n$ of quotient maps on the annulus and on the M\"obius band for $r=4$. Then, using the formula (\ref{eq:r_reg_maps_torus}), we can obtain the numbers $\tau^{(4)}_n$ of maps on the torus counted up to all its homeomorphisms. The results of these calculations are given in the Table \ref{table_last}.

\section{Toroidal 3-regular maps up to all homeomorphisms}

\begin{figure}[ht]
\centering
	\centering
    	\includegraphics[scale=0.7]{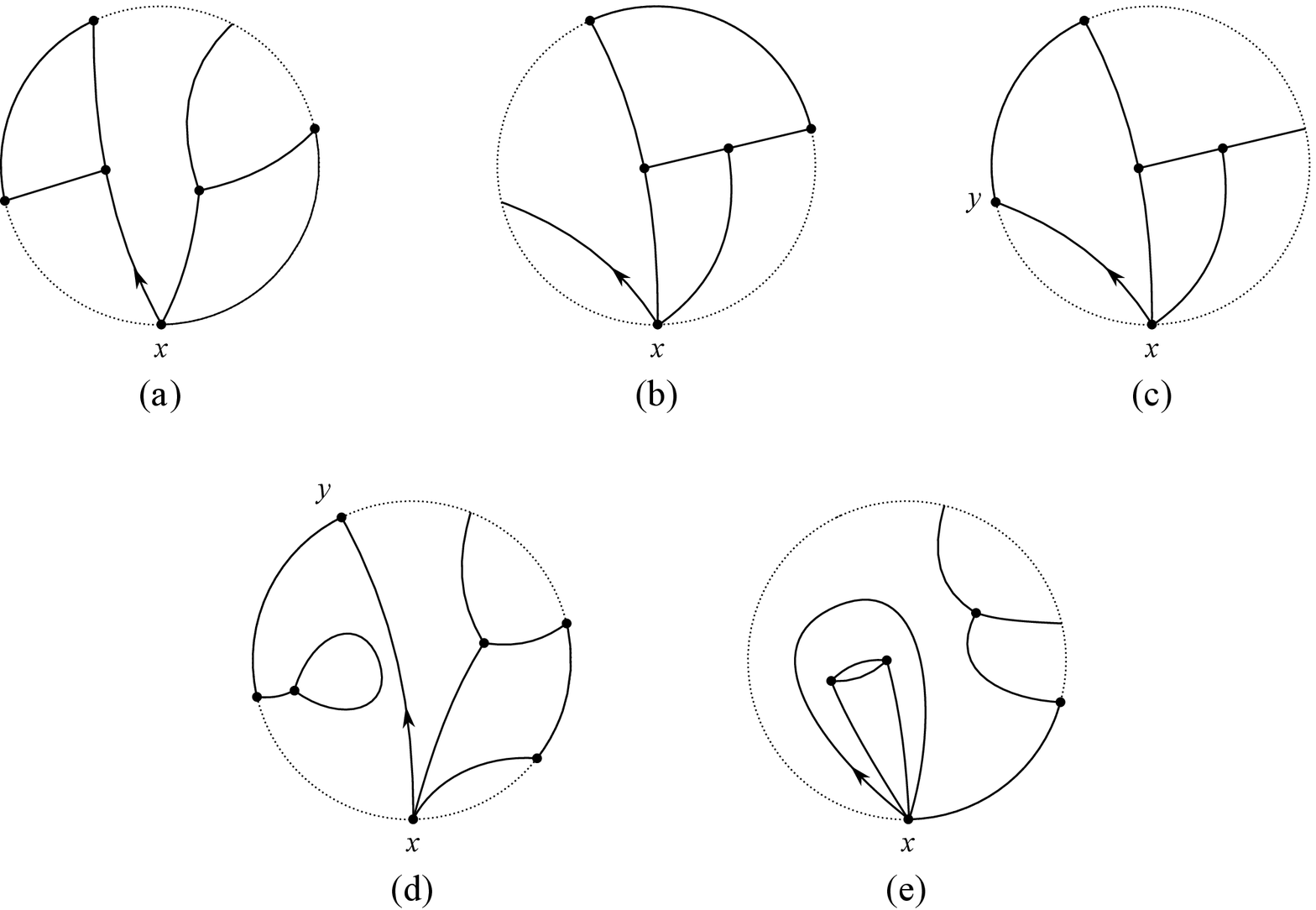}
	\caption{}
\label{fig:disk_3_regular}
\end{figure}

In the final part of the article we briefly describe analogous regarding $3$-regular maps on a torus. The differences with $4$-regular maps can be illustrated on quotient maps on the disc. As an example, we will demonstrate the derivation of a recurrence relation for the numbers $\mathfrak{d}^{\rm out}_{n-2,d-2}$ of quotient maps on the disc with the root vertex  $x$ lying on the boundary (Figure \ref{fig:disk_3_regular}).

As before, it will be enough to count such maps under the restriction that the root semi-edge is the one that precedes the boundary in a counter-clockwise ordering of semi-edges around the root vertex. The set of such quotient maps can be split into five classes (Figure \ref{fig:disk_3_regular}). The first class contains quotient maps with the root edge joining the root vertex $x$ with some other vertex lying in the interior of the disc (\ref{fig:disk_3_regular},a). The second class contains quotient maps with the root semi-edge being dangling (\ref{fig:disk_3_regular},b). The third and the fourth classes contain two types of quotient maps with the root edge joining the root vertex and a degree $2$ vertex that also lies on the boundary (Figure \ref{fig:disk_3_regular},c,d). Finally, the fifth bucket contains quotient maps with the root edge being a loop (Figure \ref{fig:disk_3_regular},e). Consequently, the recurrence relation for the numbers of such quotient maps takes the form 
\begin{equation}
\label{eq:disk_regular_3}
\mathfrak{d}^{\rm (out)}_{n,d} = \mathfrak{d}^{\rm (out)}_{n-2,d+1} +\mathfrak{d}^{\rm (out)}_{n-1,d-1} + \mathfrak{d}^{\rm (out)}_{n-3,d} + \sum\limits_{i=0}^{n-3} \mathfrak{d}^{\rm (out)}_{i,1}\cdot \mathfrak{d}^{\rm (out)}_{n-i-3,d-1} + \sum\limits_{i=0}^{n-2}\sum\limits_{j=0}^{d-2} \tilde{s}_{i,j}\cdot \mathfrak{d}^{\rm (out)}_{n-2-i,d-2-j}.
\end{equation}

Using analogous reasoning, one can prove the following recurrence relations for the numbers $\mathfrak{d}^{\rm in}_{n,d}$ that describe maps with the root vertex lying in the interior of the disc:
$$
\mathfrak{d}^{\rm (in)}_{n,d} = \mathfrak{d}^{\rm (in)}_{n-2,d+1} + \mathfrak{d}^{\rm (out)}_{n-1,d-1} + 2 \cdot \mathfrak{d}^{\rm (out)}_{n-3,d} +2\sum\limits_{i=0}^{n-2}\sum\limits_{j=0}^{d-2} \tilde{s}_{i,j}\cdot \mathfrak{d}^{\rm (in)}_{n-2-i,d-2-j}.
$$
To enumerate quotient maps on the annulus we will also need the numbers $\mathfrak{d}^{(1)}_{n,d}$ of quotient maps on the disc with the root vertex on the boundary and an additional distinguished leaf that also lies on this boundary. For these numbers, the recurrence relation takes the form
$$
\mathfrak{d}^{(1)}_{n,d} = \mathfrak{d}^{(1)}_{n-2,d+1} + \mathfrak{d}^{(1)}_{n-1,d-1} + \mathfrak{d}^{(1)}_{n-3,d} + \mathfrak{d}^{\rm (out)}_{n-2,d-1} + \sum\limits_{i=0}^{n-2} \bigl ( \mathfrak{d}^{\rm (out)}_{i,1}\cdot \mathfrak{d}^{(1)}_{n-i-2,d-1}  + \mathfrak{d}^{(1)}_{i,1}\cdot \mathfrak{d}^{\rm (out)}_{n-i-2,d-1} \bigr ) +
$$
$$
+\sum\limits_{i=0}^{n-2}\sum\limits_{j=0}^{d-2} \tilde{s}_{i,j}\cdot \mathfrak{d}^{(1)}_{n-2-i,d-2-j}.
$$

Next we provide recurrence relations for quotient maps on the annulus which can be proved by the same means. For the maps with the root vertex on the boundary the corresponding numbers $c^{\rm (out)}_{n,d}$ can be calculated by the formula
$$
c^{\rm (out)}_{n,d} = c^{\rm (out)}_{n-2,d+1} + c^{\rm (out)}_{n-1,d-1} + c^{\rm (out)}_{n-3,d} + \mathfrak{d}^{(1)}_{n-3,d-1} + \sum\limits_{i=0}^{n-3} \bigl ( \mathfrak{d}^{\rm (out)}_{i,1}\cdot c^{\rm (out)}_{n-i-3,d-1} + c^{\rm (out)}_{i,1}\cdot \mathfrak{d}^{\rm (out)}_{n-i-3,d-1} \bigr ) +
$$
$$
+\sum\limits_{i=0}^{n-2}\sum\limits_{j=0}^{d-2} \bigl (\tilde{s}_{i,j}\cdot c^{\rm (out)}_{n-2-i,d-2-j} + \mathfrak{d}^{\rm (out)}_{i,j}\cdot \mathfrak{d}^{\rm (in)}_{n-2-i,d-2-j} \bigr ).
$$
For maps on the annulus with the root vertex lying in its interior the formula is
$$
c^{\rm (in)}_{n,d} = c^{\rm (in)}_{n-2,d+1} + c^{\rm (out)}_{n-1,d-1} + 2 \cdot c^{\rm (out)}_{n-3,d} 
+ \sum\limits_{i=0}^{n-2}\sum\limits_{j=0}^{d-2}  \bigl (2\cdot \tilde{s}_{i,j}\cdot c^{\rm (in)}_{n-2-i,d-2-j} + 
 \mathfrak{d}^{\rm (in)}_{i,j}\cdot \mathfrak{d}^{\rm (in)}_{n-2-i,d-2-j}\bigr ).
$$

Finally, for maps on the M\"obius band the following formulae hold:
$$
m^{\rm (out)}_{n,d} = m^{\rm (out)}_{n-2,d+1} + m^{\rm (out)}_{n-1,d-1} + m^{\rm (out)}_{n-3,d} + 
\sum\limits_{i=0}^{n-3}\bigl (  \mathfrak{d}^{\rm (out)}_{i,1}\cdot m^{\rm (out)}_{n-i-3,d-1}+
m^{\rm (out)}_{i,1}\cdot \mathfrak{d}^{\rm (out)}_{n-i-3,d-1}
\bigr ) +  \mathfrak{d}^{(1)}_{n-3,d-1} + 
$$
$$
+ (d-1)\cdot\mathfrak{d}^{\rm (out)}_{n-2,d-2} +
\sum\limits_{i=0}^{n-2}\sum\limits_{j=0}^{d-2} \bigl (\tilde{s}_{i,j}\cdot m^{\rm (out)}_{n-2-i,d-2-j} +\mathfrak{d}^{\rm (out)}_{i,j}\cdot \tilde{p}_{n-2-i,d-2-j}\bigr ),
$$
$$
m^{\rm (in)}_{n,d} = m^{\rm (in)}_{n-2,d+1} + m^{\rm (out)}_{n-1,d-1}+ 2 \cdot m^{\rm (out)}_{n-3,d} +
(d-1) \cdot  \mathfrak{d}^{\rm (in)}_{n-2,d-2} 
+ 2\sum\limits_{i=0}^{n-2}\sum\limits_{j=0}^{d-2} \bigl(\tilde{s}_{i,j}\cdot m^{\rm (in)}_{n-2-i,d-2-j} +
 \mathfrak{d}^{\rm (in)}_{i,j}\cdot \tilde{p}_{n-2-i,d-2-j}\bigr ).
$$

The formulae (\ref{eq:cylinder_total}) and (\ref{eq:moebius_total}) for $r=3$ should be replaced by two following relations:
$$
C_n^{(3)}=c^{\rm (in)}_{n,3}+3\cdot c^{\rm (out)}_{n,2},\qquad\qquad 
M_n^{(3)}=m^{\rm (in)}_{n,3}+3\cdot m^{\rm (out)}_{n,2}.
$$
Indeed, all quotient maps on the annulus (and on the M\"obius band) can be split into two classes, depending on whether or not the root vertex lies in the interior of the surface. It is straightforward to see that maps with the root in the interior are described by the numbers $c^{\rm (in)}_{n,3}$ and $m^{\rm (in)}_{n,3}$. For maps with the root on the boundary we should use the numbers $c^{\rm (out)}_{n,2}$ and $m^{\rm (out)}_{n,2}$ and multiply them by a coefficient that describes the real number of rootings of a quotient map of the same type to remove the restriction which requires the root to be in a certain position.

Substituting these expressions into (\ref{eq:r_reg_maps_torus}) we obtain the numbers $\bar{\tau}^{(3)}_n$ of $3$-regular maps on the torus counted up to all its homeomorphisms (Table \ref{table_last}).

\section*{Conclusion}

Although for sensed maps the problem of enumerating maps by genus was solved in a generic way \cite{Mednykh_Nedela}, a similar problem for unsensed maps, being its natural generalization, seems to be little investigated. The technique developed in \cite{Liskovets_reductive_technique} for the sphere was never attempted to be used for higher genera surfaces. It seems that the technical difficulty of enumerating quotient maps on orbifolds arising when considering orientation-reversing symmetries which was noted in \cite{Liskovets_reductive_technique} would allow to obtain only quite cumbersome formulae for the general case of a genus $g$ surface. 

In this regard, the case of the torus is somewhat special. Even though we had to consider closed non-orientable surfaces (the Klein bottle) and both orientable and non-orientable surfaces with a boundary (the disc, the annulus and the M\"obius band), neither of these surfaces had branch points, and that allowed us to keep our recurrence relations reasonably complex. The results of enumeration of quotients of $r$-regular maps on these surfaces and the final formulae for unsensed toroidal $3$- and $4$-regular maps on the torus seem to be new. We have also computed the numbers of unsensed toroidal maps for $r=3 \dots 6$ and included all these results in the Table \ref{table_last}. Correctness of the first few terms in each of these sequences was checked by explicit generation of the corresponding non-isomorphic maps.

The work was supported by grant 17-01-00212 from the Russian Foundation for Basic Research. 

\begin{table}[h!]
\begin{center}
\footnotesize
\begin{tabular}{c|cccc}
\midrule
$ v $ &\phantom{00000}$\bar{\tau}^{(3)}(2v)$\phantom{00000}&\phantom{0000}$\bar{\tau}^{(4)}(v)$\phantom{0000}&\phantom{000000000000}$\bar{\tau}^{(5)}(2v)$\phantom{000000000000}&\phantom{000000000}$\bar{\tau}^{(6)}(v)$\phantom{000000000}\\ 
\midrule
1 & 1 & 1 & 13 & 3\\
2 & 5 & 4 & 3523 & 61\\
3 & 40 & 20 & 2035550 & 1936\\
4 & 450 & 133 & 1421177130 & 89986\\
5 & 6370 & 1013 & 1055597813091 & 4791784\\
6 & 104498 & 9209 & 812108624237833 & 272005507\\
7 & 1843324 & 89889 & 640086212334600319 & 15929826713\\
8 & 33778574 & 929373 & 513617627395229165708 & 951610091294\\
9 & 632053347 & 9880120 & 417872608954804473932525 & 57659992554993\\
10 & 11983323029 & 107087360 & 343735500499416537210021983 & 3532378891197016\\
\midrule
\end{tabular}
\caption{Unlabelled $r$-regular toroidal maps counted up to all homeomorphisms}
\label{table_last}
\end{center}
\end{table}

\end{document}